\documentclass[reqno]{amsart}
\usepackage{hyperref}
\usepackage{geometry}
\usepackage{enumitem}
\usepackage{amsmath}
\usepackage{amssymb}
\usepackage{amsthm}
\usepackage{amsrefs}
\usepackage{amsfonts}
\usepackage{xcolor}
\usepackage{mathtools}
\usepackage{stackengine}

\makeatletter
\@namedef{subjclassname@2020}{%
  \textup{2020} MSC}
\makeatother

\hypersetup{
	colorlinks   = true, %Colours links 
	urlcolor     = blue, %Colour for external hyperlinks
	linkcolor    = purple, %Colour of internal links
	citecolor   = blue %Colour of citations
}

\def\XXint#1#2#3{{\setbox0=\hbox{$#1{#2#3}{\int}$ }
\vcenter{\hbox{$#2#3$ }}\kern-.6\wd0}}

\makeatletter
\newcommand*{\rom}[1]{\expandafter\@slowromancap\romannumeral #1@}

\newcommand{\diam}{\text{diam}}

\newcommand{\Exact}{\text{Exact}}

\newcommand{\R}{\mathbb{R}}

\newcommand{\N}{\mathbb{N}}

\newcommand{\inv}{^{-1}}
\newcommand{\bthm}{\begin{thm}}
\newcommand{\ethm}{\end{thm}}
\newcommand{\bproof}{\begin{proof}}
\newcommand{\eproof}{\end{proof}}
\newcommand{\blem}{\begin{lem}}
\newcommand{\elem}{\end{lem}}
\newcommand{\brem}{\begin{rem}}
\newcommand{\erem}{\end{rem}}
\newcommand{\eeqn}{\end{equation}}
\newcommand{\eeqnn}{\end{equation*}}
\newcommand{\beqn}{\begin{equation}}
\newcommand{\beqnn}{\begin{equation*}}
\newcommand{\eprop}{\end{prop}}
\newcommand{\eexm}{\end{exm}}
\newcommand{\enexm}{\end{nexm}}
\newcommand{\ecor}{\end{cor}}
\newcommand{\bcor}{\begin{cor}}
\newcommand{\bexm}{\begin{exm}}
\newcommand{\bnexm}{\begin{nexm}}
\newcommand{\bprop}{\begin{prop}}
\newcommand{\bdefn}{\begin{defn}}
\newcommand{\edefn}{\end{defn}}
\newcommand{\benum}{\begin{enumerate}}
\newcommand{\eenum}{\end{enumerate}}

\title[Exact approximation order]{Exact approximation order and well-distributed sets}

\begin{document}
\theoremstyle{plain}
\newtheorem{thm}{Theorem}[section]
\newtheorem{lem}[thm]{Lemma}
\newtheorem{prop}[thm]{Proposition}
\newtheorem{cor}[thm]{Corollary}

\theoremstyle{definition}
\newtheorem{defn}[thm]{Definition}
\newtheorem{exm}[thm]{Example}
\newtheorem{nexm}[thm]{Non Example}
\newtheorem{prob}[thm]{Problem}

\theoremstyle{remark}
\newtheorem{rem}[thm]{Remark}

\author{Prasuna Bandi}
\address{\textbf{Prasuna Bandi} \\
I.H.E.S., Universit\'e Paris-Saclay, CNRS, Laboratoire Alexandre Grothendieck. 35 Route de Chartres, 91440 Bures-sur-Yvette, France}
\email{prasuna@ihes.fr}

\author{Anish Ghosh}
\address{\textbf{Anish Ghosh} \\
School of Mathematics,
Tata Institute of Fundamental Research, Mumbai, India 400005}
\email{ghosh@math.tifr.res.in}

\author{Debanjan Nandi}
\address{\textbf{Debanjan Nandi} \\
Faculty of mathematics and computer science, Weizmann Institute of Science, 234 Herzl Street, Rehovot 76100, Israel}
\email{debanjan.nandi@weizmann.ac.il}
\date{}

\thanks{A.\ G.\ gratefully acknowledges support from a grant from the Infosys foundation, a Department of Science and Technology, Government of India, Swarnajayanti fellowship and a grant from the Department of Atomic Energy, Government of India, under project $12-R\&D-TFR-5.01-0500$. D.N gratefully acknowledges support from ISF Grant 1149/18.}

\thanks{}

%\date{}

\begin{abstract} 
We prove that for any proper metric space $X$ and a function $\psi:(0,\infty)\to(0,\infty)$ from a suitable class of approximation functions, the Hausdorff dimensions of the set $W_\psi(Q)$ of all points $\psi$-well-approximable by a well-distributed subset $Q\subset X$, and the set $E_\psi(Q)$ of points that are exactly $\psi$-approximable by $Q$, coincide. This answers in a general setting, a question of Beresnevich-Dickinson-Velani in the case of approximation of reals by rationals, and answered by Bugeaud in that case using the continued-fraction expansion of reals. Our main result applies in particular to approximation by orbits of fixed points of a wide class of discrete groups of isometries acting on the boundary of hyperbolic metric spaces.
\end{abstract}

\maketitle

%\tableofcontents

\section*{Introduction}

In this article we calculate the Hausdorff dimension of ``exactly approximable" sets in a general setting. Our motivation stems from classical Diophantine approximation. Let $\psi:(0,\infty)\to(0,\infty)$ be a decreasing function such that $x \to x^2\psi(x)$ is non-increasing and consider the set of $\psi$-approximable real numbers 
\begin{equation}\nonumber
W_{\psi} := \left\{ x \in \mathbb{R}~:~\left|x - \frac{p}{q}\right| < \psi(q) \text{ for infinitely many } \frac{p}{q}\right\}.
\end{equation}
Khintchine \cite{Khintchine} proved that $W_{\psi}$ has Lebesgue measure zero if and only if 

\begin{equation}\label{sumcong}
\sum_{x \geq 1} x\psi(x) < \infty.
\end{equation}

Now let $W_{\tau} := W(x \to x^{\tau})$. Then Jarn\'{i}k \cite{Jar1} and independently Besicovitch \cite{Besicovitch} proved that $$\dim (W_{\tau}) = 2/\tau.$$ Under the assumption that $\psi(x) = o(x^{-2})$, Jarn\'{i}k in \cite{Jar2} further proved that the set
\begin{equation}\nonumber
E_\psi := W_{\psi}\backslash \bigcup_{0<c<1} W_{c\psi}
\end{equation}
 is non-empty. This result was refined by G\"{u}ting \cite{Gut}, who proved that for $\tau \geq 2$, $$\dim E_\tau = 2/\tau = \dim (W_{\tau}).$$
Throughout this paper, dim will mean Hausdorff dimension.

 The set $\Exact(\psi)$ above is the set of real numbers $x$ such that
 \begin{equation}\label{exact1}
 |x - p/q| \leq \psi(q) \text{ infinitely often}
 \end{equation}
 and
 \begin{equation}\label{exact2}
 |x - p/q| \geq c\psi(q) \text{ for any } c < 1 \text{ and any } q \geq q_0(c, x).
 \end{equation}
 
More generally, one can consider the set $\Exact(\psi, \phi)$, which is defined by replacing $c\psi$ by $\phi$ in (\ref{exact2}) above.

 Beresnevich, Dickinson and Velani \cite{BDV} (see also \cites{BDVmem, DV}) studied $\Exact(\psi, \phi)$ under certain conditions; their methods use in particular that $\psi/\varphi\to \infty$. 
 
 In \cite{Bu1} Bugeaud proved the following result, resolving a conjecture of Beresnevich, Dickinson and Velani from \cite{BDV}.
 
 \bthm \label{Bugeaud}
  Let $\psi:(0,\infty)\to(0,\infty)$ be such that $x \to x^2\psi(x)$ is non-increasing. Assume that $\sum_{x \geq 1} x\psi(x)$ converges and let $\lambda$ denote the lower order at infinity of $1/\psi$. Then 
$$ \dim (E_\psi) = 2/\lambda = \dim (W_{\psi}).  $$
 \ethm
 
Here the lower order at infinity $\lambda(g)$ of a function $g:(0,\infty)\rightarrow (0,\infty)$ is defined by 
\begin{equation*}
    \lambda(g)=\liminf_{x\rightarrow \infty}\frac{\log g(x)}{\log x}.
\end{equation*}

There has been extensive further work on this subject by
 Bugeaud (\cites{Bu2, Bu3}), Moreira \cite{Moreira}, and Bugeaud and Moreira \cite{BuMo} who have established comprehensive metric results for $\Exact(\psi)$, including the situation when $\sum_{x \geq 1} x\psi(x)$ diverges. In addition to their intrinsic interest, the sets $\Exact(\psi)$ are also connected to other aspects of Diophantine approximation, for instance the Lagrange spectrum \cite{BuMo}. Besides the classical Diophantine approximation  problem of reals by rationals, there are several other situations  where a theory of metric Diophantine approximation exists. A rich class of examples arises from the action of discrete groups of isometries on boundaries of hyperbolic spaces. 
This raises the natural question.

 \begin{center}
 \textit{Is it possible to develop a general framework to investigate $\Exact(\psi)$?}\\
\end{center} 

 In the case of rationals in the real line, Bugeaud uses classical results from the theory of continued fractions to construct a Cantor set whose points are $\psi$-exactly approximable and whose Hausdorff dimension has the required lower bound in terms of $\psi$, to prove Theorem \ref{Bugeaud}. We show that this is a metric phenomenon, that is, such a Cantor set of exactly approximable points exists, given any countable set that is suitably well separated and distributed, in appropriate metric spaces.  
 
Note that although exact-approximation is a completely metric notion the continued-fraction expansion of reals and related approximation properties with respect to the euclidean metric in itself is not sufficient for a lower bound on the Hausdorff dimension of exactly approximable points. A key requirement as is usual in Hausdorff lower bound estimates, is the existence of measures with suitable metric properties. In the real line one has the Lebesgue measure and it is not immediate what are the appropriate metric conditions for such a measure in more general settings. 

We now introduce the general setting. Let $(X,d)$ be a metric space. Let $Q\subset X$ be a countable set and $R:Q\to \mathbb{R}^+$ a function. For such a pair $(Q,R)$, and a function $\psi: \R \rightarrow \R^+$, define the $\psi$-well approximable set with respect to $(Q,R)$ to be
$$ W_{\psi}(Q,R):=\{x\in X : d(x,\xi)<\psi(R(\xi)) \text { for infinitely many }\xi  \in Q \}, $$
and define the exact set to be
$$E_{\psi}(Q,R):= W_{\psi}(Q,R)\setminus \bigcup_{0<c<1}W_{c\psi}(Q,R). $$
Define the lower order at zero of $\psi$ to be $$\underline{\lambda}(\psi):=\liminf_{x\rightarrow 0}\frac{\log \psi(x)}{\log x}.$$
We now state the main result of the paper, postponing the relevant definitions, which are a bit technical, to the next section (Definitions \ref{regular_measure} and \ref{WD}).

\bthm \label{main0}
Let $(X,d,\mu)$ be an $(\alpha,\beta)$-regular metric space. Let $Q$ be a WDS subset of $X$ with radius function $R$. Let $\psi: (0,\infty) \rightarrow (0,\infty)$ be a non-decreasing function such that $\sum_{\xi \in Q}\left(\frac{\psi(R(\xi))}{R(\xi)}\right)^{\alpha}$ converges. 
Then 
\begin{equation}\label{main}
\dim E_{\psi}(Q,R) = \dim(W_{\psi}(Q, R)) = \frac{\alpha}{\underline{\lambda}(\psi)}. 
\end{equation}
\ethm

\subsection{Remarks}
\begin{enumerate}
\item The substance of the above theorem is the inequality
$$ \dim E_{\psi}(Q,R) \geq \frac{\alpha}{\underline{\lambda}(\psi)}.$$
The upper bound follows from relatively straight forwards covering arguments, while $\dim(W_{\psi}(Q, R))$ has been computed in greater generality in \cite{GhoNan}. In fact the study of the dimension of $W_{\psi}(Q, R)$ has an older vintage. We refer the reader to \cite{GhoNan} for the history of the problem. 
\item We apply our results to Diophantine approximation of limits of suitable discrete groups acting properly discontinuously on hyperbolic spaces by orbits of parabolic and loxodromic fixed points (see \S \ref{example}), a subject which was systematically first treated in Patterson's landmark paper \cite{Patterson}, and has seen considerable activity recently \cites{BGSV, GhoNan, HV} and the references therein.
\item For the Jarn\'{i}k-Besicovitch type problems, one usually considers regular or well distributed systems to study the approximation problem. However, the exact set up needs an even stronger assumption (Definition \ref{WD}).
\item In an upcoming paper we consider the sets of exact order in the situation of Diophantine approximation by arbitrary open sets as in the earlier work of \cite{GhoNan} of the second and third named authors. This has applications to questions on the asymptotic properties of geodesics in negatively curved spaces.
\item Note that the dimension lower bound in Theorem \ref{main0} is the same as in Bugeaud's Theorem \ref{Bugeaud}, (the function $\psi$ considered in Bugeaud \cite{Bu1} is a function of the height of a rational number in $[0,1]$, whereas our $\psi$ is a function of a radius function $R$ which in the real case is the reciprocal of height squared).
\end{enumerate}

\section{Preliminaries}

\bdefn[Regular measure]\label{regular_measure} Let $(X,d)$ be a metric space. A Borel measure $\mu$ on $X$ is called $(\alpha,\beta)$-regular, where $\alpha,\beta >0$, if there exists $C>0$ such that for each ball $B=B(\xi,r)\subset X$ and each annulus $A(\xi,r,R):=\overline{B}(\xi,R)\setminus B(\xi, r)\subset X$, where $0\leq r\leq R\leq \diam(X)$, the following hold:
\benum
\item (Ahlfors regularity). $$r^\alpha/C\,\leq \,\mu(B)\,\leq\, C\cdot r^\alpha,$$
\vspace{0.2cm}
\item (Annular regularity). $$\frac{1}{C}\cdot r^\alpha\cdot \left(\frac{R-r}{r}\right)^{\beta}\,\leq\, \mu(A(\xi,r,R)).$$
\eenum
\edefn

We note below that essentially, in the Ahlfors regular case, the condition of regularity (for some $0<\beta\leq \alpha$) is equivalent to annuli being non-empty. 
\blem\label{annular}
Suppose $(X,d,\mu)$ is $\alpha$-Ahlfors regular. Then for each $\xi\in X$ and $0\leq r<R\leq \diam(X)$, the annulus $A(\xi,r,R)\neq \emptyset$ iff $(X,d,\mu)$ is $(\alpha,\alpha)$-regular.
\elem
\bproof
The claim follows by Ahlfors-regularity and $\emptyset\neq \overline{B}(\xi,\frac{3R}{8}+\frac{5r}{8})\backslash B(\xi,\frac{R}{8}+\frac{7r}{8})\subset \overline{B}(\xi, R)\backslash B(\xi,r)$.
\eproof
\bdefn[WDS subsets]\label{WD} Let $(X,d,\mu)$ be an $\alpha$-Ahlfors regular metric measure space. A countable subset $Q\subset X$ is WDS with radius-function $R:Q\to \R_+$, if there exists a constant $C\geq 1$ such that the following hold.
\benum
\item (\textbf{well-separated}). For $\xi,\zeta\in Q$, 
\begin{equation*}
    d(\xi,\zeta)\,\geq\, \frac{1}{C}\cdot \min\,\{R(\xi),R(\zeta)\}.
\end{equation*}

\vspace{0.2cm}
\item (\textbf{well-distributed}). For each ball $B=B(\eta,r)\subset X$, there exists a number $k_B=k_B(\eta,r)\in \N$, such that for each $k\geq k_B$, there is a subset $Q_B(k)\subset Q$, where
\vspace{0.2cm}
\benum
\item for all $\xi\in Q_B$, $1/(C\cdot k) \,\leq\, R(\xi)\,\leq\, C/k,$
\vspace{0.2cm}
\item $\bigcup_{\xi\in Q_B} B(\xi, R(\xi))\subset B,$
\vspace{0.2cm}
\item $\#Q_B \geq \frac{1}{C}\cdot k^\alpha \cdot \mu(B)$.
\vspace{0.2cm}
\eenum
\eenum
\edefn

\blem\label{opensets}Let $(X,d,\mu)$ be a proper $\alpha$-Ahlfors regular metric measure space. If $(Q,R)$ is WDS then requirements (1) and (2) of WDS hold for any bounded open sets with quantitative and uniform bounds on constants, that is, given $U\subset X$ open, there is a constant $0<k_U\in\N$ such that for $k\geq k_U$, properties (2a), (2b) and (2c) are satisfied.
\elem
\bproof
We assume without loss of generality that $U$ is a connected non-empty open sets. 
Write for $\epsilon>0$ $$U_\epsilon=\{\xi\in U\,|\, d(\xi,X\setminus U)\geq \epsilon\}.$$ 
Choose an $\epsilon_0>0$ such that $\mu(U_{\epsilon_0})\geq \mu(U)/2$. By a $5r$-covering lemma, pick $r_0=r_0(\epsilon_0)>0$, and a collection of disjoint balls $\{B(\xi_i,r_0)\}_i$ of radius $r_0$, with centers in $U_{\epsilon_0}$, such that $$U_{\epsilon_0}\,\subset\,\bigcup_i B(\xi_i,5\cdot r_0)\,\subset\, U.$$ The lemma follows by Ahlfors regularity and applying the definition of well-distributedness to the balls $B(\xi_i,r_0)$ individually. 
\eproof

\subsection{Hausdorff measure}
\bdefn
Let $(X,d)$ be a metric space. Suppose $F\subseteq X$ and let $\rho >0$. A $\rho$-cover for $F$ is any countable collection of balls $(B_{i})_{i\in \N}$ with $r(B_{i})\leq \rho$ for all $i\in \N$ and $F\subseteq \bigcup_{i\in \N}B_{i}$. For $s\geq 0$, define
\begin{equation*}
    \mathcal{H}^{s}_{\rho}(F)=\inf \left\{\sum_{i=1}^{\infty}r(B_{i})^{s} : \{B_{i}\} \text{ is a } \rho \text{-cover of }F \right\}
\end{equation*}
The Hausdorff $s$-measure of $F$ is then defined as
\begin{equation*}
    \mathcal{H}^{s}(F):=\lim_{\rho\rightarrow 0}\mathcal{H}^{s}_{\rho}(F)
\end{equation*}
and Hausdorff dimension of $F$ is defined as 
\begin{equation*}
    \dim F :=\inf \{s:\mathcal{H}^{s}(F)=0\}=\sup \{s:\mathcal{H}^{s}(F)=\infty\}.
\end{equation*}
\edefn
\blem(Mass Distribution Principle) \\
Let $m$ be a probability measure supported on a subset $F$ of $X$. Suppose there are positive constants $c$ and $r_{0}$ such that
\begin{equation*}
    m(B)\leq c r^{s}
\end{equation*}
for any ball $B$ with radius $r\leq r_{0}$. Then $\dim F \geq s$.
\elem
Next we state a general version of a well-known application of the mass-distribution principle which we will use. The proof is straightforward generalization of Example 4.6 in \cite[pg. 64]{Falcon}.

\blem \label{prop:1} Let $X$ be a proper $\alpha$-regular metric space. Let $\{K_{l}=\bigsqcup_{i=1}^{t_{l}}C_{i}^{l}\}_{l\geq 0}$ be a nested decreasing sequence of finite disjoint unions of compact sets. Suppose the following holds:
 \begin{enumerate}
     \item For $l\geq 1$, each compact set of $K_{l}$ contains at least $m_{l}\geq 2$ compact sets of $K_{l}$
     \item the maximal diameter of the compact sets in $K_{l}$ tends to $0$ as $l\rightarrow \infty$.
     \item 
    For $l\geq 1$, $d(C_{i}^{l},C_{j}^{l})\geq \delta_{l} \text{ where }0<\delta_{l+1}<\delta_{l} \text{ and } 1\leq i\neq j \leq t_{l}$
 \end{enumerate}
Then 
 \begin{equation*}
     \dim \bigcap _{l=0}^{\infty}K_{l}\geq \alpha \liminf_{l\rightarrow \infty}\frac{\log (m_{1}\cdots m_{l-1})}{-\log \left(m_{l}\delta_{l}^{\alpha}\right)}.
 \end{equation*}
 \elem

\subsection{Hyperbolic metric spaces}\label{example} As we see below a large class of examples of well-distributed systems and regular measures come from looking at actions of discrete subgroups of isometry groups of hyperbolic metric spaces and this section is devoted to discussing them.
A geodesic metric space $(X,\rho)$ is $\delta$-Gromov hyperbolic for $\delta\geq 0$, if any geodesic triangle is $\delta$-thin, that is any side is contained in the $\delta$-neighbourhood of the union of the two other sides, see \cite{BH} for basic notions regarding these spaces. We will denote by $\partial X$ the Gromov boundary of $X$. When $X$ is proper, $\partial X$ is a compact metric space (with topology equivalent to the topology as the ideal boundary), equipped for each $x\in X$ with visual metrics
$$d^x_\epsilon(\xi,\eta) \approx e^{-\epsilon\cdot(\xi|\eta)_x},$$ for $\xi,\eta\in \partial X$, where $\approx$ denotes bilipschitz with absolute constants, and
$$(\xi|\eta)_x=\underset{\substack{w\to\xi\\z\to\eta}}{\lim}\,\, \frac{1}{2}\cdot(\rho(w,x)+\rho(z,x)-\rho(w,z)),$$ is the Gromov product extended to the boundary (see \cite{BH}). 
Usually $x$ will be fixed beforehand and we will just write $d$ for the corresponding visual metric.

Suppose now that $\Gamma$ is a discrete subgroup of the group of isometries of a proper geodesic hyperbolic space $(X,\rho)$ that acts on $X$ properly discontinuously and cocompactly. In this case there is a unique (upto constant multiples) $\Gamma$-equivariant quasiconformal density of dimension $v_\Gamma$ (see \cite{Coo}), on the boundary, which is ergodic for the induced action (by quasisymmetries) of $\Gamma$ on $\partial X$, which is known as the Patterson-Sullivan density $\mu_\Gamma$, where $v_\Gamma$ is the critical exponent for the action of $\Gamma$ defined
$$v_\Gamma=\limsup_{r\to 0}\,\frac{\log(\#\{g\in\Gamma: \rho(x,gx)\leq r\})}{r},$$ the definition being independent of $x$,
and by a $\Gamma$-invariant quasiconformal density of dimension $\alpha$ in $\partial X$, one means a family $\{\nu_x\}_{x\in\Gamma}$ of measures in $\partial X$ such that the following hold:
\benum
\item ($\Gamma$-equivariance) For all $x\in X$, $g\in \Gamma$, $g_\ast \nu_x= \nu_{gx}$, where $g_\ast\nu_x$ is the pushforward $\nu_x\circ g\inv$.
\item (Quasiconformality) For all $x,y\in Y$, $\nu_y$-a.e $\xi\in \partial X$, $$
\frac{d\nu_x}{d\nu_y}(\xi)\approx e^{-\alpha\cdot\beta_\xi(x,y)};$$
\eenum
here $$\beta_\xi(x,y)=\sup_{\{z_n\}}\limsup_{z_n\to\xi} \,(\rho(x,z_n)-\rho(y,z_n)),$$ is the Busemann function. By Sullivan's shadow lemma (see \cite{Sullivan}, \cite{Quint}) the measure $(X,d,\mu_\Gamma)$ is $v_\Gamma$-Ahlfors-regular in this case.

If $\Gamma$ is of divergence type, that is, $$\sum_{g\in\Gamma}e^{-v_\Gamma\cdot \rho(x,gx)}=\infty,$$
again the Patterson-Sullivan density is the unique $v_\Gamma$-dimensional conformal density in the limit set $\Lambda_\Gamma$ and is ergodic for the action of $\Gamma$ (see \cite{Roblin}). This is the case when the set of conical limit points of $\Gamma$ in $\partial X$ have full $\mu_\Gamma$-measure. The density $\mu_\Gamma$ is not in general Ahlfors-regular with respect to the visual metric $d$ restricted to its support $\Lambda_\Gamma$. If the space has negative curvature (CAT) bounded above and below, and the action of $\Gamma$ is additionally assumed to be geometrically finite, then $\mu$ is doubling (see \cite{SV}, \cite{Sch}). 

We note below a sufficient condition on $(X,\Gamma)$ which implies the regularity of the quasiconformal density (using Lemma \ref{annular}).

\bcor\label{annularex}
Let $(X,\rho)$ be a proper, geodesic hyperbolic metric space and assume that $\Gamma$ is a discrete group acting non-elementarily and properly discontinuously by isometries on $X$. Let $\mu$ be a $\Gamma$-quasiconformal density of dimension $\alpha>0$, which is $\alpha/\epsilon_0$-Ahlfors-Regular in $\Lambda_\Gamma$ equipped with the visual metric $d_\infty(\cdot,\cdot)\approx e^{-\epsilon_0(\cdot|\cdot)_{x_0}}$, where $\epsilon_0=\epsilon_0(X)>0$ is a constant. Assume that for all $\xi\in \partial X$ and $0\leq r<R<\diam(\Lambda_\Gamma)$, $(\overline{B}(\xi, R)\backslash B(\xi,r))\cap \Lambda_\Gamma$ is non-empty. Then $\mu_{x_0}$ is an $(\alpha/\epsilon_0,\alpha/\epsilon_0)$-regular measure.
\ecor

\brem[Examples] \benum We discuss here some examples.
\item Instances where Corollary \ref{annularex} applies include lattices on rank one symmetric spaces. The Ahlfors-regularity in this case follows from Sullivan's shadow lemma (and ergodicity of the action in case of non-uniform lattices) \cite{Sullivan}. More generally, geometrically finite actions with Ahflors-regular Patterson-Sullivan densities and connected limit sets satisfy the hypothesis. This includes for example suitable quasi-Fuchsian actions. The orbits of loxodromic and parabolic limit points in these cases are WDS with suitable radius function $R$. Thus the hypothesis of Theorem \ref{main0} are satisfied. Namely, with $\rho$ as the metric of the symmetric space, if $\xi$ is the parabolic fixed point of the maximal parabolic subgroup $\Gamma_\xi$ of a lattice $\Gamma$, then $Q=\Gamma(\xi)$ and $R(g\xi)=C_x e^{-\rho(x,g'x)}$, where $C_x$ is a constant depending only on $x$, and $g'$ is a coset representative of $g\Gamma_\xi$, which minimizes the distance of $x$ to $g\Gamma_\xi(x)$ (see \cite{Sullivan}). In this case the measure is $(D_X,D_X)$-regular, where $D_X$ is the asymptotic logarithmic growth of balls in the space. If $\xi$ is a loxodromic fixed point of the maximal subgroup $L_\xi$, then again $R(g\xi)=C_xe^{-\rho(x,g'x)}$, where $g'$ is the coset representative of $gL_\xi$ which minimises distance to $x$.

\item In general, if $\Gamma$ has a geometrically finite action on a proper hyperbolic metric space, with virtually nilpotent parabolic subgroups, then there always exist $\Gamma$-conformal densities which are Ahlfors-regular for suitable visual metrics on the limit set $\Lambda_\Gamma$, and $\Gamma$-orbits of parabolic and loxodromic fixed points are WDS; the Ahlfors-regularity in this case follows from a Sullivan-type estimate for shadows of balls in a suitable hyperbolic space with $\Lambda_\Gamma$ as visual boundary, the WDS property comes from geometrical finiteness of the action (see \cite{Nandi}). In this case as in the previous one, the radius function also is defined in terms of nearest distances to the relevant cosets, with respect to suitable hyperbolic metrics.

\item The second property required in Definition \ref{regular_measure}, of course does not follow automatically for the Patterson-Sullivan measure, even when $X$ admits a properly discontinuous, cocompact action. Consider the Cayley graph of the free group in two generators (which is a tree and a CAT($-1$) space). The Patterson-Sullivan measure for the action of this group on this Cayley graph does not satisfy the second requirement of regularity. It is easy to see that the conclusion of Theorem \ref{main0} also fails in this case, with $Q$ as the orbit of any loxodromic fixed point $\xi$, and radius function as mentioned above (namely $e^{-\|g'\|}$ for the point $g\xi$, where $\|g'\|$ is the number of letters in the minimal word representing $g'$ and $\|g'\|$ is minimal such that $g'\xi=g\xi$). 
\eenum
\erem

\section{Proof of main result}
\bproof
Fix a ball $B=\overline{B}(x,r)$ in $X$ such that $\mu (B)<\infty$. For $l\in \N$, let $c_{l}=1-1/2^{l}$.\\
For $\xi$ in $Q$ and $l\in \N$, define
$$A_{l}(\xi):= A(\xi, c_{l}\psi(R(\xi)),\psi(R(\xi))).$$
Using induction, we will construct a rapidly increasing sequence of integers $(N_{l})_{l=1}^{\infty}$ and a decreasing sequence of subsets $B=K_{0}\supseteq K_{1}\supseteq \cdots \supseteq K_{n}\supseteq \cdots $ of $B$  such that 
\begin{enumerate}[label=\textbf{P.\arabic*}]
\item \label{p1} for each $l\in \N$,
\begin{equation*}
    K_{l}=\bigsqcup_{i=1}^{t_{l}}A_{l}(\xi_{i}^{l})
\end{equation*}
for some $\xi_{i}^{l} \in Q$, $t_l\in\N$, with $\frac{1}{CN_{l}} \leq R(\xi_{i}^{l})\leq \frac{C}{N_{l}} $.
\item \label{p2} For each $l\in \N $ and $1\leq i\neq j\leq t_{l}$,
\begin{equation*}
    d(\xi_{i}^{l},\xi_{j}^{l})\geq 4\max \{\psi(R(\xi_{i}^{l})): 1\leq i \leq t_{l}\}.
\end{equation*}
\item  \label{p3} The number of annuli of $K_{l+1}$ that lie in a fixed annulus of $K_{l}$ is at least
\begin{equation}
    m_{l+1}:= \begin{cases}
    n_{l+1}N_{l+1}^{\alpha} \psi\left(\frac{1}{CN_{l}}\right)^{\alpha} & \text{if $l\geq 1$}\\
    \frac{1}{C}N_{1}^{\alpha}\mu(B) & \text{if $l=0$}
    \end{cases}
\end{equation}
where $n_{l+1}=\frac{1}{2C^{2}}c_{l}^{\alpha-\beta}(1-c_{l})^{\beta}$.\\
\item \label{p4} For $l\geq 1$ and $x\in K_{l}$, the inequality $d(x,\xi)<\psi(R(\xi))$ has at least $l$ solutions $\xi \in Q$ with $ R(\xi)\geq \frac{1}{CN_{l}}$. 
\item \label{p5} For $l\geq 2$ and $x\in K_{l}$, the inequality $d(x,\xi)<c_{l-1} \psi(R(\xi))$ has no solution $\xi \in Q$ with $ \frac{1}{CN_{l}} \leq R(\xi)\leq \frac{C}{N_{l-1}} $. 
\end{enumerate}

Since $(Q,R)$ is well separated, for all $\xi,\eta \in Q$ we have
\begin{equation}
    \tag{WS}
    d(\xi,\eta)\geq \frac{1}{C}\min\{R(\xi),R(\eta)\}.
    \label{eqn:WS}
\end{equation}

\noindent\textbf{Constructing $K_{1}$}:\\
Denote by $k_{B}$ the constant arising by applying well distributedness of $(Q,R)$ to the ball $B$. Let $\varepsilon_{0}$ be small enough such that for all $\varepsilon\leq \varepsilon_{0}$
\begin{equation}\label{eqn:14}
    \frac{\psi(\varepsilon)}{\varepsilon}<\frac{1}{4C^{3}}.
\end{equation}
Such an $\varepsilon_{0}$ exists since $\psi(x)/x \rightarrow 0$ as $x\rightarrow 0$. Let
\begin{equation}\label{eqn:21}
    N_{1}\geq \max\{k_{B},C/\varepsilon_{0}\}.
\end{equation}
 By well distributedness, there exists elements $\xi_{1}^{1},\xi_{2}^{1},\ldots, \xi_{t_{1}}^{1}$ (depending on $N_{1}$) of $Q$ such that
 \begin{equation}\label{eqn:15}
     1/(C N_{1}) \,\leq\, R(\xi_{i}^{1})\,\leq\, C/N_{1} \text{ for all } 1\leq i\leq t_{1},
 \end{equation}
 \begin{equation}\label{eqn:16}
     \bigcup_{i=1}^{t_{1}}B(\xi_{i}^{1},R(\xi_{i}^{1}))\subseteq B,
 \end{equation}
 and 
 \begin{equation}\label{eqn:17}
     t_{1} \geq \frac{1}{C}\cdot N_{1}^\alpha \cdot \mu(B).
 \end{equation}
Since $$R(\xi_{i}^{1})\stackrel{\eqref{eqn:15}}{\leq}C/N_{1}\stackrel{\eqref{eqn:21}}{\leq}\varepsilon_{0}$$ by \eqref{eqn:14} we get 
\begin{equation}\label{eqn:30}
    \psi(R(\xi_{i}^{1}))\leq \frac{R(\xi_{i}^{1})}{4C^{3}}\stackrel{\eqref{eqn:15}}{\leq}\frac{1}{4C^{2}N_{1}} \;\text{ for all}\; 1\leq i\leq t_{1}.
\end{equation}
For $i\neq j$, $A_{1}(\xi_{i}^{1})\cap A_{1}(\xi_{j}^{1})=\emptyset$. Indeed, if $\gamma \in A_{1}(\xi_{i}^{1})\cap A_{1}(\xi_{j}^{1})$, then
\begin{align*}
    \frac{1}{C^{2}N_{1}}\stackrel{\eqref{eqn:15}}{\leq} \frac{1}{C}\min\{R(\xi_{i}^{1}),R(\xi_{j}^{1})\}& \stackrel{\eqref{eqn:WS}}{\leq} d(\xi_{i}^{1},\xi_{j}^{1})\\
    &\leq d(\xi_{i}^{1},\gamma)+d(\xi_{j}^{1},\gamma)\\
    &\leq \psi(R(\xi_{i}^{1}))+\psi(R(\xi_{j}^{1}))\stackrel{\eqref{eqn:30}}{\leq} \frac{1}{2C^{2}N_{1}}
\end{align*}
which is a contradiction.
 Now define
$$K_{1}:= \bigsqcup_{i=1}^{t_{1}} A_{1}(\xi_{i}^{1}).$$
Then 
\begin{equation*}
    K_{1}\stackrel{\eqref{eqn:30}}{\subseteq} \bigcup_{i=1}^{t_{1}}B(\xi_{i}^{1},R(\xi_{i}^{1}))\stackrel{\eqref{eqn:16}}{\subseteq} B.
\end{equation*}
 For $1\leq i\neq j \leq t_{1}$
 \begin{align*}
   4 \max\{\psi(R(\xi_{i}^{1})):1\leq i\leq t_{1}\}&\stackrel{\eqref{eqn:30}}{\leq}\frac{1}{C^{2}N_{1}} \\ &\stackrel{\eqref{eqn:15}}{\leq} \frac{1}{C}\min\{R(\xi_{i}^{1}),R(\xi_{j}^{1})\} \stackrel{\eqref{eqn:WS}}{\leq} d(\xi_{i}^{1},\xi_{j}^{1})
 \end{align*}
 thus proving \ref{p2}. By definition of $K_{1}$ and \eqref{eqn:15}, \ref{p1} and \ref{p4} are satisfied. Finally, note that \eqref{eqn:17} implies \ref{p3}.\\ \\
\textbf{Constructing $K_{l+1}$}:\\
Suppose $K_{l}$ is constructed satisfying \ref{p1}, \ref{p4} and $N_{l}$ is chosen large enough such that 
\begin{equation}\label{eqn:22}
     \sum_{\{\xi \in Q \; : \; R(\xi)<\frac{C}{N_{l}}\}}R(\xi)^{-\alpha}\psi(R(\xi))^{\alpha}<\frac{(1-c_{l})^{\beta}}{(8C^{3})^{2\alpha +3}c_{l}^{\beta-\alpha}}. 
\end{equation}
Such an $N_{l}$ exists since $\sum_{\xi \in Q}\left(\frac{\psi(R(\xi))}{R(\xi)}\right)^{\alpha}$ converges.\\
Denote by $\mathcal{A}_{l}:=\{\xi_{i}^{l}:1\leq i \leq t_{l}\}$.
Fix an annulus from $K_{l}$, say $A_{l}(\xi)$ for some $\xi \in \mathcal{A}_{l}$  with 
\begin{equation}\label{eqn:3}
    1/(C N_{l}) \leq R(\xi)\leq C/N_{l}.
\end{equation}
\textbf{Step 1}:
For $a\in (0,\infty)$ and $\xi \in \mathcal{A}_{l}$, define 
$$ S_{a}(\xi):=\{\eta \in Q : R(\eta)= a \text{ and } B(\eta, \psi(a))\cap A_{l}(\xi)\neq \emptyset \}.$$
Claim: For $a <\frac{C}{N_{l}}$, we have 
\begin{equation}\label{eqn:31}
|S_{a}(\xi)|\leq (8C^{3})^{1+\alpha}\dfrac{c_{l}^{\beta-\alpha}}{(1-c_{l})^{\beta}}a^{-\alpha}\mu(A_{l}(\xi)).
\end{equation}
For $\eta_{1},\eta_{2} \in S_{a}(\xi), \eta_{1}\neq \eta_{2}$, by \eqref{eqn:WS} we get  $d(\eta_{1},\eta_{2})\geq \frac{a}{C}$. Hence the balls $\{B(\eta,\frac{a}{4C^{3}})\}_{\eta \in S_{a}(\xi)}$ are disjoint. For $a< \frac{C}{N_{l}}$, we will show that 
\begin{equation}\label{eqn:7}
    \bigsqcup_{\eta \in S_{a}(\xi)}B\left(\eta,\frac{a}{4C^{3}}\right)\subset B\left(\xi,2 \psi(R(\xi))\right).
\end{equation}
Let $\gamma \in B\left(\eta,\frac{a}{4C^{3}}\right)$ for some $\eta \in S_{a}(\xi)$. Since $\eta \in S_{a}(\xi)$, $B(\eta, \psi(a))\cap A_{l}(\xi)\neq \emptyset$. Let $\gamma' \in B(\eta, \psi(a))\cap A_{l}(\xi)$. Since $$a<C/N_{l}\leq C/N_{1}\stackrel{\eqref{eqn:21}}{\leq}\varepsilon_{0}$$by \eqref{eqn:14} we get
   \begin{equation}\label{eqn:4}
       \psi (a)<\frac{a}{4C^{3}}.
   \end{equation}
  Also, by \eqref{eqn:3}
\begin{equation}\label{eqn:2}
  R(\eta)=a< \frac{C}{N_{l}}\leq C^{2}R(\xi). 
\end{equation}
Since $\eta, \xi \in Q$ by \eqref{eqn:WS}, we have
\begin{align}\label{eqn:8}
    \frac{a}{C^{3}}=\frac{1}{C^{3}}R(\eta)&\stackrel{\eqref{eqn:2}}{\leq} d(\xi,\eta)\\ \nonumber
    &\leq d(\xi,\gamma')+d(\gamma',\eta)\\ \nonumber
    &\leq \psi(R(\xi)) +\psi(a)\\ \nonumber
    &\stackrel{\eqref{eqn:4}}{\leq} \psi(R(\xi)) + \frac{a}{4C^{3}}
\end{align}
which implies
\begin{equation}\label{eqn:5}
    \frac{a}{4C^{3}}\leq \frac{1}{3} \psi\left(R(\xi)\right).
\end{equation}
 Hence,
  \begin{align}\label{eqn:6}
      d(\xi,\gamma)&\leq \nonumber d(\xi,\gamma')+d(\gamma',\eta)+d(\eta,\gamma)\\ \nonumber
      & \leq \psi\left(R(\xi)\right)+\psi(a)+\frac{a}{4C^{3}}\\
      &\stackrel{\eqref{eqn:4},\eqref{eqn:5}}{\leq} 2 \psi\left(R(\xi)\right).
  \end{align}
This proves (\ref{eqn:7}).
Using $(\alpha,\beta)$-regularity, it follows from (\ref{eqn:7}) that for $a<\frac{C}{N_{l}}$
\begin{equation}
    |S_{a}(\xi)|\frac{a^{\alpha}}{C(4C^{3})^{\alpha}}\leq \sum _{\eta \in S_{a}(\xi)}\mu \left(B\left(\eta,\frac{a}{4C^{3}}\right)\right) \leq \mu(B(\xi,2 \psi(R(\xi))))\leq C2^{\alpha} \psi(R(\xi))^{\alpha}\leq C^{2}2^{\alpha}\frac{c_{l}^{\beta-\alpha}}{(1-c_{l})^{\beta}}\mu(A_{l}(\xi))
\end{equation}
which implies that
\begin{equation*}
    |S_{a}(\xi)|\leq (8C^{3})^{1+\alpha}\frac{c_{l}^{\beta-\alpha}}{(1-c_{l})^{\beta}}a^{-\alpha}\mu(A_{l}(\xi))
\end{equation*}
thus proving the claim.\\
\textbf{Step 2}: For each $\xi \in \mathcal{A}_{l}$, denote by $k_{l}(\xi):=k_{A_{l}(\xi)}$ the constant arising by applying well distributedness of $(Q,R)$ to $A_{l}(\xi)$ (Lemma \ref{opensets}).
Denote by 
\begin{equation*}
  k_{l}:=  \max\{k_{l}(\xi):\xi \in \mathcal{A}_{l} \}.
\end{equation*}
Choose $\varepsilon_{l}$ small enough such that for all $\varepsilon\leq \varepsilon_{l}$ 
\begin{equation}\label{eqn:32}
  \frac{\psi(\varepsilon)}{\varepsilon}<\frac{1-c_{l}}{C^{3}}  
\end{equation}
and let  
\begin{equation}\label{eqn:27}
    N_{l+1}\geq \max \{k_{l},C/\varepsilon_{l},C^{2}N_{l}\}.
\end{equation}
Then there exists $Q_{l}(\xi):=Q_{A_{l}(\xi)}(N_{l+1})=\{ \gamma_{1},\gamma_{2},\cdots,\gamma_{t_{l}(\xi)}\} \subset Q$ such that 
\begin{equation}\label{eqn:25}
    \frac{1}{CN_{l+1}} \leq R(\gamma_{i})\leq \frac{C}{N_{l+1}} \text{ for all } 1\leq i\leq t_{l}(\xi).
\end{equation}
\begin{equation*}
    \bigcup_{i=1}^{t_{l}(\xi)}B(\gamma_{i},R(\gamma_{i}))\subseteq A_{l}(\xi)
\end{equation*}
and 
\begin{equation}\label{eqn:35}
   t_{l}(\xi) \geq \frac{1}{C} N_{l+1}^\alpha \mu(A_{l}(\xi)). 
\end{equation}
Define 
$$ J(\xi):=\left\{ x\in A_{l}(\xi): d(x,\eta)\geq\psi(R(\eta)) \text{ for all }\eta \in Q \text{ with } R(\eta)<\frac{C}{N_{l}}\right\} $$ 
and
 $$ R':=\frac{1}{3C}\min \{R(\gamma): \gamma \in Q_{l}(\xi)\cap (A_{l}(\xi)\setminus J(\xi))\}.$$
Then $\{B(\gamma,R')\}_{\gamma \in Q_{l}(\xi)\cap (A_{l}(\xi)\setminus J(\xi))}$ are disjoint.
We will now check that 
\begin{equation}\label{eqn:10}
    \bigsqcup_{\gamma \in Q_{l}(\xi)\cap (A_{l}(\xi)\setminus J(\xi))}B(\gamma,R')\subset \bigcup_{\{\eta \in Q \; :\; R(\eta)<C/N_{l}\}}\bigcup_{\eta' \in S_{R(\eta)}(\xi)} B(\eta',2\psi(R(\eta'))).
\end{equation}
Let $\gamma' \in B(\gamma,R')$ for some $\gamma \in A_{l}(\xi)\setminus J(\xi)$. Since \begin{equation*}
    A_{l}(\xi)\setminus J(\xi)=\bigcup_{\{\eta \in Q \; :\; R(\eta)<C/N_{l}\}}B(\eta,\psi(R(\eta)))\cap A_{l}(\xi)
     \subseteq \bigcup_{\{\eta \in Q \; :\; R(\eta)<C/N_{l}\}}\bigcup_{\eta' \in S_{R(\eta)}(\xi)} B(\eta',\psi(R(\eta'))),
\end{equation*}
there exists $\eta \in Q$ with $R(\eta)<C/N_{l}$ and $\eta' \in S_{R(\eta)}(\xi)$ such that $\gamma \in B(\eta',\psi(R(\eta')))$. Then
$$\frac{1}{C}\min\{R(\gamma),R(\eta)\}\leq d(\gamma,\eta')<\psi(R(\eta'))=\psi(R(\eta))\stackrel{\eqref{eqn:14}\eqref{eqn:21}}{<} \frac{R(\eta)}{4C}$$ 
which implies
$$\min \{R(\gamma),R(\eta)\}=R(\gamma)$$
and hence
$$ R'\leq \frac{1}{3C}R(\gamma)\leq \frac{d(\gamma,\eta')}{3}<\frac{\psi(R(\eta'))}{3}.  $$
Therefore,
$$d(\gamma',\eta')\leq d(\gamma',\gamma)+d(\gamma,\eta')\leq R'+\psi(R(\eta')) \leq 2\psi(R(\eta'))  $$
thus proving \eqref{eqn:10}. Using Ahlfors regularity and \eqref{eqn:10}, we get
\begin{align*}
    |Q_{l}(\xi)\cap (A_{l}(\xi)\setminus J(\xi))|\frac{1}{3^{\alpha}C^{2\alpha+1}N_{l+1}^{\alpha}}&\leq \sum_{\gamma \in Q_{l}(\xi)\cap (A_{l}(\xi)\setminus J(\xi))} \mu(B(\gamma,R')) \\
    &\leq \sum_{\{\eta \in Q \; :\; R(\eta)<C/N_{l}\} } \sum_{\eta' \in S_{R(\eta)}(\xi)}\mu (B(\eta',2\psi(R(\eta'))))\\
        &\stackrel{\eqref{eqn:31}}{\leq} C2^{\alpha}(8C^{3})^{1+\alpha} \frac{c_{l}^{\beta-\alpha}}{(1-c_{l})^{\beta}}\mu(A_{l}(\xi)) \sum_{\{\eta \in Q \; : \; R(\eta)<\frac{C}{N_{l}}\}}R(\eta)^{-\alpha}\psi(R(\eta))^{\alpha} 
\end{align*} 
which implies 
\begin{equation*}
    |Q_{l}(\xi)\cap (A_{l}(\xi)\setminus J(\xi))| \leq (8C^{3})^{2\alpha +2} N_{l+1}^{\alpha}\mu(A_{l}(\xi))\frac{c_{l}^{\beta-\alpha}}{(1-c_{l})^{\beta}}\sum_{\{\eta \in Q \; : \; R(\eta)<\frac{C}{N_{l}}\}}R(\eta)^{-\alpha}\psi(R(\eta))^{\alpha}  \stackrel{\eqref{eqn:22}}{\leq} \frac{1}{2C}N_{l+1}^{\alpha}\mu(A_{l}(\xi)).
\end{equation*}
 Let $Q'_{l}(\xi)=Q_{l}(\xi)\cap J(\xi)$. 
Then
\begin{equation*}
   |Q'_{l}(\xi)|= t_{l}(\xi)-|Q_{l}(\xi)\cap (A_{l}(\xi)\setminus J(\xi))|\stackrel{\eqref{eqn:35}}{\geq} \frac{1}{2C}N_{l+1}^{\alpha}\mu(A_{l}(\xi)).
\end{equation*}
Hence, using $(\alpha,\beta)$-regularity and $\psi$ is non-decreasing, we get that for any $\xi \in \mathcal{A}_{l}$,
\begin{equation}\label{eqn:23}
    |Q'_{l}(\xi)|\geq n_{l+1}N_{l+1}^{\alpha} \psi\left(\frac{1}{CN_{l}}\right)^{\alpha}
\end{equation}
where $n_{l+1}=\frac{1}{2C^{2}}c_{l}^{\alpha-\beta}(1-c_{l})^{\beta}$.\\
Using a similar argument as in construction of $K_{1}$, one can show that 
$$A_{l+1}(\gamma_{1})\cap A_{l+1}(\gamma_{2})=\emptyset \;\; \text{ for }\;\; \gamma_{1}\neq \gamma_{2}\in Q_{l}'(\xi) $$
and 
\begin{equation}\label{eqn:24}
    A_{l+1}(\gamma)\subseteq A_{l}(\xi) \text{ for any }\gamma \in Q_{l}'(\xi).
\end{equation}
Denote by 
\begin{equation*}
    \{\xi_{i}^{l+1}:1\leq i\leq t_{l+1}\}:=\bigsqcup_{\xi \in \mathcal{A}_{l}}Q_{l}'(\xi)
\end{equation*}
where $t_{l+1}:=\sum_{\xi \in \mathcal{A}_{l}}|Q'_{l}(\xi)|$ and
define
\begin{equation*}
    K_{l+1}:=\bigsqcup_{i=1}^{t_{l+1}}A_{l+1}(\xi_{i}^{l+1}).
\end{equation*}
By $\eqref{eqn:24}$, we have $K_{l+1}\subset K_{l}$. \eqref{eqn:25} implies that \ref{p1} holds for $K_{l+1}$. It is easy to see that \ref{p2} holds using a similar argument as in the construction of $K_{1}$ and \ref{p3} holds by \eqref{eqn:23}. \\ \\
\textbf{Step 3}: 
In this step we show that \ref{p4} and \ref{p5} hold for $K_{l+1}$.\\

Let $\zeta \in K_{l+1}$. Then there exists $\xi \in \mathcal{A}_{l}$ and $\gamma\in Q'_{l}(\xi)$ such that $\zeta \in A_{l+1}(\gamma).$ 
Since $K_{l+1}\subseteq K_{l}$, by  induction hypothesis the equation 
\begin{equation}\label{eqn:26}
    d(\zeta,\eta)\leq \psi(R(\eta))
\end{equation}
 has at least $l$ solutions $\eta\in Q$ with $\frac{1}{CN_{l}}\leq R(\eta)$. Since $\zeta \in A_{l+1}(\gamma)$, $\gamma$ is also a solution of the equation \eqref{eqn:26} and it cannot be one of the above $l$ solutions since $R(\gamma)\leq \frac{C}{N_{l+1}}\stackrel{\eqref{eqn:27}}{<}\frac{1}{CN_{l}}\leq R(\eta)$. Hence the equation \eqref{eqn:26} has at least $l+1$ solutions $\eta\in Q$ with $\frac{1}{CN_{l+1}}\leq R(\eta)$. This proves \ref{p4}.\\
 
Suppose there exists $\eta \in Q$ with $\frac{1}{CN_{l+1}}\leq R(\eta)< \frac{C}{N_{l}}$ and $d(\zeta, \eta)<c_{l}\psi(R(\eta))$. Since $\gamma \in J(\xi)$ and $R(\eta)< \frac{C}{N_{l}}$, we have $d(\gamma, \eta)\geq \psi(R(\eta))$.
Since $\gamma \in Q_{l}(\xi)$, $R(\gamma)\stackrel{\eqref{eqn:25}}{\leq}C/N_{l+1}\stackrel{\eqref{eqn:27}}{\leq}\varepsilon_{l}$ and hence by \eqref{eqn:32}
\begin{equation}\label{eqn:33}
    \psi(R(\gamma))<\frac{(1-c_{l})}{C^{3}}R(\gamma)
\end{equation}
By \eqref{eqn:25}, $R(\gamma)\leq \frac{C}{N_{l+1}}\leq C^{2}R(\eta)$. Hence,
\begin{align}\label{eqn:11}
    \frac{1}{C^{3}}R(\gamma)\leq \frac{1}{C}\min \{R(\eta),R(\gamma)\}\nonumber
    \leq d(\eta,\gamma) \nonumber
    &\leq d(\eta,\zeta)+d(\zeta,\gamma)\\ \nonumber
    &\leq c_{l}\psi(R(\eta))+\psi(R(\gamma))\\ 
    &\stackrel{\eqref{eqn:33}}{<} c_{l}\psi(R(\eta))+\frac{(1-c_{l})}{C^{3}}R(\gamma)
\end{align}
which implies 
\begin{equation} \label{eqn:12}
    R(\gamma)< C^{3}\psi(R(\eta)).
\end{equation}
By (\ref{eqn:11}) and (\ref{eqn:12}), we have
\begin{equation*}
    d(\eta,\gamma)< c_{l}\psi(R(\eta))+\frac{(1-c_{l})}{C^{3}}R(\gamma) < \psi(R(\eta))
\end{equation*}
contradicting that $\gamma \in J(\xi)$. This proves \ref{p5}.\\
Thus we have constructed a decreasing sequence of subsets $(K_{l})_{l=1}^{\infty}$ of $B$ satisfying all the five properties stated above.
Now, define
\begin{equation*}
    \mathbb{K}:=\bigcap_{l=1}^{\infty}K_{l}
\end{equation*}
\ref{p4} and \ref{p5} imply that $\mathbb{K}\subseteq E_{\psi}(Q,R)$. Hence, providing lower bound on $\dim \mathbb{K}$ gives a lower bound on $\dim E_{\psi}(Q,R)$.\\ \\
\textbf{Calculating Hausdorff dimension}:
We use lemma \ref{prop:1} to get a lower bound on $\dim \mathbb{K}$.\\
For each $l\in \N$, let $\delta_{l}=\frac{1}{3}\min \{d(\xi_{i}^{l},\xi_{j}^{l})-\psi(R(\xi_{i}^{l}))-\psi(R(\xi_{i}^{l})) : 1\leq i\neq j\leq t_{l}\}$. Then it is easy to see that $d(A_{l}(\xi_{i}^{l}),A_{l}(\xi_{j}^{l}))\geq \delta_{l}$. Since $$R(\xi_{i}^{l})\stackrel{\ref{p1}}{\leq}C/N_{l}\leq C/N_{1}\stackrel{\eqref{eqn:21}}{\leq}\varepsilon_{0}$$ by \eqref{eqn:14} we get 
\begin{equation}\label{eqn:34}
    \psi(R(\xi_{i}^{l}))\leq \frac{R(\xi_{i}^{l})}{4C^{3}}\stackrel{\ref{p1}}{\leq}\frac{1}{4C^{2}N_{1}} \;\text{ for all}\; 1\leq i\leq t_{l}.
\end{equation}
Hence
\begin{align*}
    \delta_{l}=&\frac{1}{3}\min \{d(\xi_{i}^{l},\xi_{j}^{l})-\psi(R(\xi_{i}^{l}))-\psi(R(\xi_{i}^{l})) : 1\leq i\neq j\leq t_{l}\}\\
    &\stackrel{\eqref{eqn:WS}\eqref{p1}\eqref{eqn:34}}{\geq}\frac{1}{3}\left(\frac{1}{C^{2}N_{l}}-\frac{1}{2C^{2}N_{1}}\right)\stackrel{}{=} \frac{1}{6C^{2}N_{l}}.
\end{align*}
Therefore 
\begin{equation*}
    m_{l+1}\delta_{l+1}^{\alpha} \stackrel{\ref{p3}}{\geq} \frac{1}{(6C^{2})^{\alpha}}n_{l+1}\psi\left(\frac{1}{CN_{l}}\right)^{\alpha}.
\end{equation*}
Fix $0<\varepsilon <1$. 
By the definition of $\underline{\lambda}(\psi)$, there exists a strictly decreasing sequence of real numbers $(X_{l})_{l\geq 1}$ tending to $0$ such that
\begin{equation}\label{eqn:29}
    \log \psi\left(X_{l}\right)\geq (\underline{\lambda}(\psi) +\varepsilon)\log X_{l} \;\;\text{ for all }l\geq 1.
\end{equation}
We may assume that $(\frac{1}{CX_{l}})_{l\geq 1}$ is a strictly increasing sequence of positive integers. We now assume additionally that $(N_{l})_{l\geq 1}$ is a subsequence of $(\frac{1}{CX_{l}})$ and that
\begin{equation}\label{eqn:40}
    N_{l} \geq \max \left\{\left(C n_{l}^{-1/\alpha}\psi\left(\frac{1}{CN_{l-1}}\right)^{-1}\right)^{2/\varepsilon},\left((6C^{2})^{\alpha}n_{l+1}^{-1}\right)^{\frac{2(1-\varepsilon)}{\alpha \varepsilon}} \right\}.
\end{equation}
Then
\begin{align*}
    \dim \mathbb{K} &\geq \alpha\liminf_{l\rightarrow \infty}\frac{\log m_{l}}{-\log (m_{l+1}\delta_{l+1}^{\alpha})} \\
    & \stackrel{\eqref{eqn:40}}{\geq} \alpha \liminf_{l\rightarrow \infty}\frac{(\alpha -\alpha\varepsilon/2)\log CN_{l}}{-\log ((6C^{2})^{-\alpha}n_{l+1})-\alpha \log \psi(\frac{1}{CN_{l}})}\\
    &\stackrel{\eqref{eqn:29},\eqref{eqn:40}}{\geq} \alpha\liminf_{l\rightarrow \infty}\frac{\alpha(1-\varepsilon)\log CN_{l}}{\alpha (\underline{\lambda}(\psi)+\varepsilon)\log CN_{l}} =\frac{\alpha(1-\varepsilon)}{\underline{\lambda}(\psi)+\varepsilon}.
\end{align*}
Thus $\dim \mathbb{K}\geq \alpha/\underline{\lambda}(\psi)$.

\eproof

\brem For Diophantine approximation on the real line, one often restricts attention to the unit interval, exploiting the translation invariance of the problem. We remark that this is a general phenomenon.

Let $Y$ be a locally-compact, Hausdorff space. Let $\Gamma$ be a countable group acting freely, properly on $Y$ by homeomorphisms. Then there is a fundamental domain $F_\Gamma=\overline{U}_\Gamma$ for the action of $\Gamma$, where $U_\Gamma\subset Y$ is open. \\
We indicate the construction: Consider $[y]\in \Gamma\setminus Y$. Consider the projection $p:Y\to \Gamma\backslash Y$ and choose $y\in p\inv([y])$. Pick an open neighbourhood $U_y$ of $y$ such that $gU_y\cap U_y=\emptyset$ for $g\neq e$ (action is free and proper) and $\overline{U}_y$ is compact. Then note that $p(U_y)$ is open and $U_y\xrightarrow{p} p(U_y)$ is a homeomorphism. For each $[y]\in\Gamma\backslash Y$, choose a neighbourhood $p(U_y)$ containing it, as above. By compactness, $\Gamma\backslash X=\cup_{i=1}^{n_\Gamma} p(U_{y_i}).$

Set for $1\leq i< n_\Gamma$, $V_{y_i}:=U_{y_i}\setminus \cup_{i=j+1}^{n_\Gamma}\cup_{g\in\Gamma}g\overline{U}_{y_j}.$ By compactness of $U_{y_i}$, $V_{y_i}$ is open. Set $V_{y_{n_\Gamma}}=U_{y_{n_\Gamma}}$. Note that $\Gamma\backslash Y=p(\cup_{i=1}^{n_\Gamma}\overline{U}_{y_i})=p(\cup_{i=1}^{n_\Gamma}\overline{V}_{y_i})$ and $p:\cup_{i=1}^{n_\Gamma}\overline{V}_{y_i}\to\Gamma\backslash Y$ is injective. Write $U_\Gamma=\cup_{i=1}^{n_\Gamma} V_{x_i}$. Then Int$(\overline{U}_\Gamma\cap g\overline{U}_\Gamma)=\emptyset$ for $g\neq e$ and there is a homeomorphism $\hat{p}:\Gamma\backslash\overline{U}_\Gamma\to\Gamma\backslash Y$. It follows that $\overline{U}_\Gamma$ is a fundamental domain in $Y$ for the action of $\Gamma$. Now we check that if $\Gamma\backslash Y$ is connected, then $\overline{U}_\Gamma$ may be chosen to be connected. Choose a metric $\hat{d}$ which generates the topology of $\Gamma\backslash Y$. Then note that $\cup_{g\in\Gamma}g\overline{V}_{x_i}=\overline{\cup_{g\in\Gamma}gV_{x_i}}$. As $\Gamma\backslash Y$ is connected, we have for some $i\geq 2$, $(\cup_g g\overline{V}_{x_1})\cap (\cup_g g\overline{V}_{x_j})\neq\emptyset$. Set $\hat{V}_{x_2}=g_2V_{x_2}$, where $\overline{V}_{x_1}\cap g_2\overline{V}_{x_i}\neq\emptyset$. Write $\hat{V}_{x_1}=V_{x_1}$. By induction we get for each $1< j\leq n_\Gamma$, $\hat{V}_{x_j}$, a translate of some $V_{x_k}$ such that $(\cup_{l=1}^{j-1}\overline{\hat{V}}_{x_l}\cap \overline{\hat{V}}_{x_j}\neq\emptyset$ and all $\{\hat{V}_{x_i}\}_{1\leq i\leq n_\Gamma}$ are mutually disjoint. The set $\overline{U}_\Gamma$, where $U_\Gamma=\cup_{i=1}^{n_\Gamma}\hat{V}_{x_i}$ is the required set.

Then note that when $Y=\partial X\setminus \infty$, where $\infty$ is the fixed point of a maximal parabolic subgroup $\Gamma_\infty$ of a discrete group $\Gamma$, acting geometrically finitely on the proper geodesic hyperbolic space $X$, then one can consider consider a fundamental domain $F_\Gamma=\overline{U}_\Gamma$ for the action of $\Gamma_\infty$ on $\partial X\setminus\infty$ and $(gF_\Gamma,d_\infty)$, and $(F_\Gamma,d_\infty)$ are bilipschitz equivalent for any $g\in \Gamma_\infty$ and visual metric $d_\infty$. 
\erem

%%%%%%%%%%%%%%%%%%%%%%%%%%%%%%%%%%%%%%%%%%%%%%

\end{document}